\documentclass[preprint,12pt]{elsarticle}

\graphicspath{ {./figures/} }
\usepackage{subfig}
\usepackage{amssymb}

\begin{document}

\begin{frontmatter}

\title{Novel Transition to fully absorbing state without long-range spatial 
	order in Directed Percolation class}

	\author[1]{Sumit S. Pakhare}
	\ead{sumitspakhare@gmail.com}
	\author[2]{Prashant M. Gade\corref{cor1}}
	\ead{prashant.m.gade@gmail.com}
	\cortext[cor1]{Corresponding Author}
	\address[1]{Department of Physics, Kamla Nehru Mahavidyalaya, Nagpur}
	\address[2]{Department of Physics, RTM Nagpur University, Nagpur}
\begin{abstract}
We study coupled Gauss maps in one dimension and observe a transition to 
band periodic state with 
2 bands. This is a periodic state with period-2 in a coarse-grained 
	sense.
This state does not show any long-range order in 
space. We compute two different order parameters to quantify the transition
	a) Flipping rate $F(t)$ which measures departures from period-2 and
	b) Persistence $P(t)$ which quantifies the loss of memory of initial 
	conditions. At the critical point, $F(t)$ shows a power-law decay 
	with exponent 0.158 which is close to 1-D directed percolation (DP)
	transition. The persistence
exponent at the critical point is found to be 1.51
which matches with several models 
in 1-D DP class. We also study  the finite-size scaling and off-critical 
	scaling to estimate other  exponents $z$ and $\nu_{\parallel}$.
We observe excellent scaling for both $F(t)$ as well as $P(t)$ and
the exponents obtained are clearly in DP class.
We believe that DP transition could be observed in systems where 
activity goes to zero even if the spatial profile could be inhomogeneous
and lacking any long-range order. 
\end{abstract}


\end{frontmatter}


\section{Introduction}
\label{}
Nonequilibrium phase transitions are found in a variety of situations.
In nature, we observe
dynamic  phases ranging from the synchronized flashing of fireflies,
spirals or rings in a chemical reaction, synchronized chirping of 
crickets,
Turing patterns etc.
The studies in this
fascinating field have only started. The most studied transition
from the viewpoint of phase transitions is the
transition to absorbing state.
They are further  classified in classes such as
directed percolation, compact directed percolation, parity conserving class,
voter class, Manna universality class, etc \cite{henkel2008non} .
Directed Ising universality class has been observed
in systems such as Grassberger's model A and B \cite{grassberger1984new,grassberger1989some},
branching and annihilating random walks with two offsprings
\cite{takayasu1992extinction,sudbury1990branching,jensen1993conservation,jensen1994critical}, 
interacting monomer-dimer model \cite{kim1994critical} and nonequilibrium
kinetic Ising model \cite{menyhard1994one,menyhard1995non}. In one dimension, 
voter universality class in voter model is equivalent to 
compact directed  percolation.
It can be mapped to an equilibrium model \cite{henkel2008non}.
There has been a long-standing debate if pair contact process
with diffusion (PCPD) is a new universality class 
 \cite{hinrichsen2001pair,noh2004universality}.
Directed percolation remains most studied and
most observed universality class in this context.
Even for PCPD, increasing
evidence points to the possibility that for long enough simulations on
large enough systems, it will be in directed percolation universality 
class \cite{mahajan2018stretched,matte2016persistence}.

If we consider models with continuous variable values  such
as coupled map lattices, there have been fewer studies.
Transition in
logistic map with delay (when mapped on a pseudo spatiotemporal system)
is found to be in directed Ising
class \cite{mahajan2013dynamic,jabeen2006spatiotemporal}. There
are systems which share the transition to equilibrium systems.
The transition to an antiferromagnetic
state in coupled logistic maps is found to be in Glauber-Ising
class \cite{gade2013universal}.
These systems are studied in higher dimension as well.
In two dimensions, the possibility of transitions in equilibrium,  as well
as nonequilibrium class has been studied. Work by Miller
and Huse demonstrating the possibility of transition in Ising class
in coupled map lattice in two dimensions attracted a lot of
attention \cite{miller1993macroscopic}. It was found that the nature of update
matters and
it can change universality class.
Coupled map lattices with some specific maps been shown
to be in $q=3$ Potts class \cite{salazar2005critical}.
Chat{\'e} and
Mannevile studied the transition to a laminar phase
in coupled piecewise linear discontinuous maps and showed that
the transition is in directed percolation universality
class \cite{chate1988spatio}.
Transition to spatiotemporal intermittency in coupled circle maps
is also found to be in directed percolation universality
class \cite{jabeen2006spatiotemporal,menon2003persistence}. 
Chat{\'e} and Mannevile studied
the transition to spatiotemporal intermittency in 2-D coupled maps and
showed that the continuous transition is in directed percolation
universality class \cite{chate1988continuous}.
In many of these cases, change in nature of update can change universality class.
Transitions  in DP class are often transitions to a synchronized state
and the state has long-range order.
Even for transition is to a chaotic synchronized state, where
infinite absorbing states are possible, there is an obvious long-range 
order \cite{PhysRevE.73.036212}.
The above transitions are marked by clean order parameters such as
the number of active sites or number of domain walls.

Janssen and Grassberger conjectured that
the transition generically belongs to directed percolation universality class
if a) transition is to a unique absorbing state from a fluctuating
active phase b) characterized
by non-negative one-component order parameter c) couplings are short-range
and d) the system has no special attributes
such as unconventional symmetries, conservation laws or quenched randomness
\cite{janssen1999levy}. 
In particle systems, directed percolation was observed 
in Domany-Kinzel automata 
\cite{hinrichsen1998numerical}, threshold
transfer processes and  Ziff-Gulari-Barshad model \cite{albano2001numerical}. 

We would demonstrate a possibility of
DP transition in a
system  which does not follow Janssen-Grassberger conjecture
not only in the sense that it does not have a unique absorbing state.
The state is not synchronized or periodic in space
and 
does not have a long-range order. There is no long-range 
spatial order
even if we coarse-grain the variables, {\it{i.e.}} 
if the variable values are divided 
into classes depending upon their values.
Furthermore, this transition does not belong to 
the damage-spreading class.
This is a transition to a frozen state in coarse-grained period-2.
This should open up the possibility of observing DP transition
experimentally in a system which approaches a   
state which is periodic in time in a coarse-grained sense.
Most studies in nonequilibrium phase transitions are in particle systems
or systems in which variables take discrete values. The
model studied in this paper 
has variables which take continuous values. Usually,
the order parameter is obtained from the spatial profile at
a given time instance such as the 
variance of the profile. For a particle system, it could be
the  number of isolated particles  or active particles at a
given time. In our case, the 
appropriate order parameter 
is obtained by observing the difference
between the spatial profile at a given time step and previous time steps.

\section{Model}
The system consists of diffusively coupled Gauss maps. Gauss map
is given as,
\[ f(x)=\exp(-\nu x^{2})+\beta,\hspace{20pt} x \in \mathbb{R}\]
Where $\nu$ and $\beta$ are the parameters.
The above function is a Gaussian with variance proportional to $1/\nu$ 
shifted by $\beta$. The value of the function tends
to $\beta$ as $x \rightarrow \infty$.
We fixed the value of $\nu = 7.5$ while $\beta$ is our control parameter.
Unlike logistic or tent map, this is not a map on the interval. 
The nature of the function for a few different values of $\beta$ is shown 
in Fig. 1. The number of fixed points changes from 3 to 1 as we 
increase $\beta$. For large values of $\beta$, the largest fixed point
is stable.  The bifurcation diagram for a single Gauss map is shown
in Fig. 2.

\begin{figure}
\centering
\includegraphics[width=80mm]{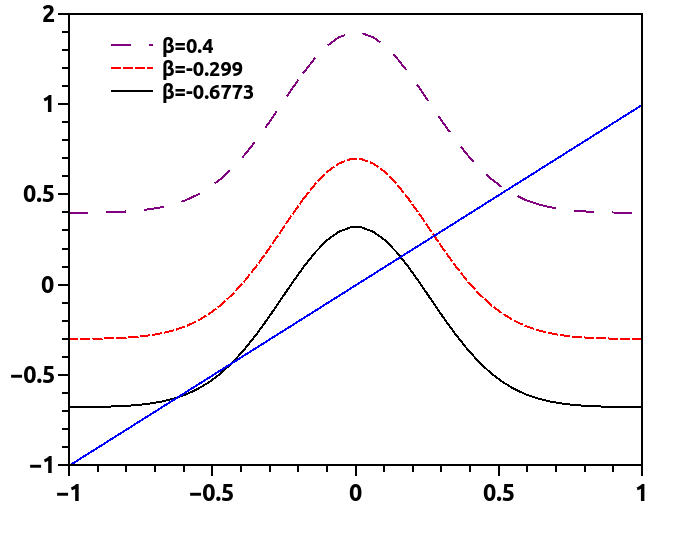}%
\caption{\label{fig:1} Gauss Map for various values of $\beta$.}
\end{figure}

\begin{figure}
\centering
\includegraphics[width=80mm]{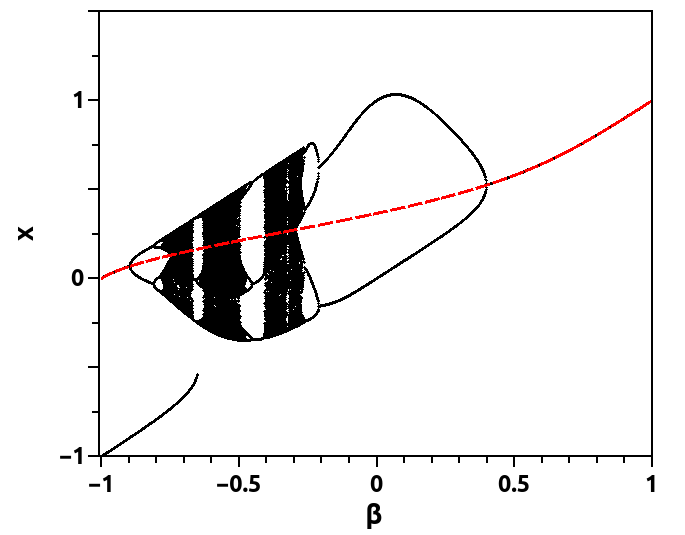}
\caption{\label{fig:2} Bifurcation Diagram of Single Gauss Map. The fixed point is
indicated by a red dashed line.  }
\end{figure}

We couple Gauss maps diffusively as follows,
\[ x_{i}(t+1)=(1-\epsilon )f(x_{i}(t))+
{\frac{\epsilon}{2}}[f(x_{i+1}(t))+f(x_{i-1}(t))]\]
where $x_{i}(t) $ is the variable value associated with the site i at time t.
We assume the periodic boundary conditions. 
We fix the coupling, $\epsilon=0.4$ and vary $\beta$. The bifurcation diagram 
for $N=100$ is shown in Fig. 3.

\begin{figure}
\centering
\includegraphics[width=80mm]{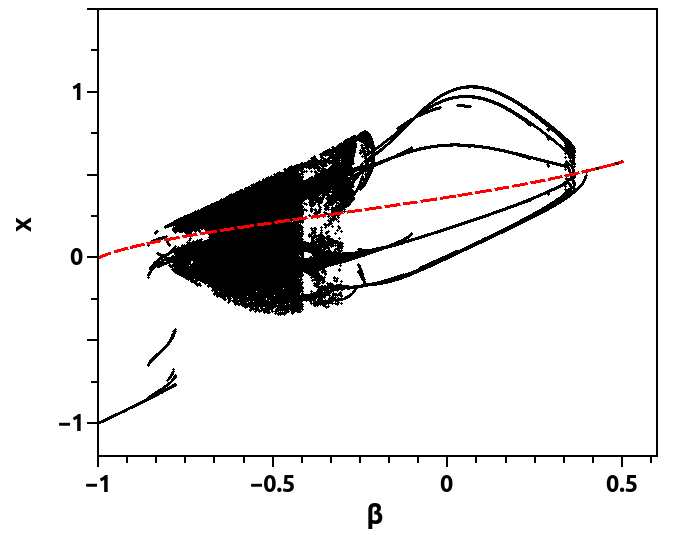}%
\caption{\label{fig:3} Bifurcation diagram of Coupled Gauss Maps. The fixed
        point is indicated by a red dashed line.}
\end{figure}

We have also shown the largest fixed point of the map as a reference. 
As mentioned above, this map
 $f(x)=\exp(-7.5 x^{2})+\beta$ has a stable fixed point for large values
of $\beta$. Simple stability
analysis indicates that
the coupled map lattice also has a stable fixed point for large values
of $\beta$ \cite{gade1993spatially}. The transition to synchronization has 
been studied extensively in many works. We will investigate another
transition in detail in this work.
A clear  two-band structure is seen for smaller  values of $\beta$. 
Not only there is a two-band structure but the system is also frozen in this state
for smaller values of $\beta$. There is a coarse-grained 2-periodicity. We
have shown the spatial profile at two different time steps for $\beta=-0.65$
and $\beta=-0.69$ in Fig. 4 and 5. We have also shown the largest fixed point for respective values of $\beta$ as reference.

\begin{figure}
\centering
\includegraphics[width=80mm]{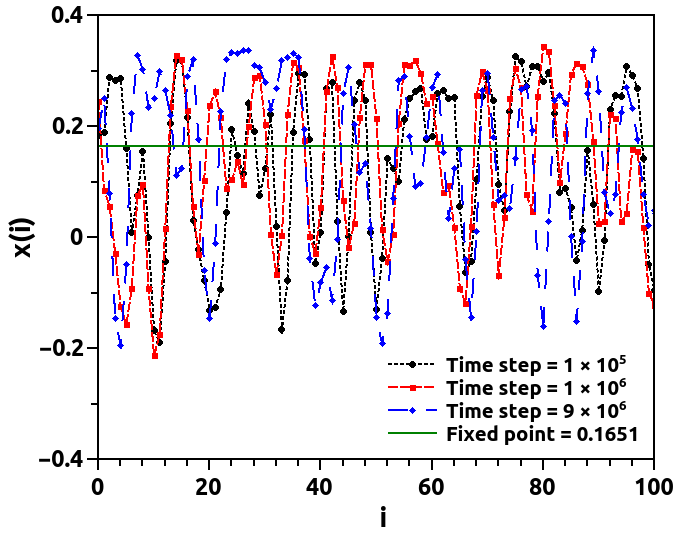}%
\caption{\label{fig:4} Non-Frozen spatial
pattern for $\epsilon=0.4$ and $\beta=-0.65$. 
	The plot is done after an even number of time steps.
	The fixed point is also shown as a reference.}
\end{figure}

\begin{figure}
\centering
\includegraphics[width=80mm]{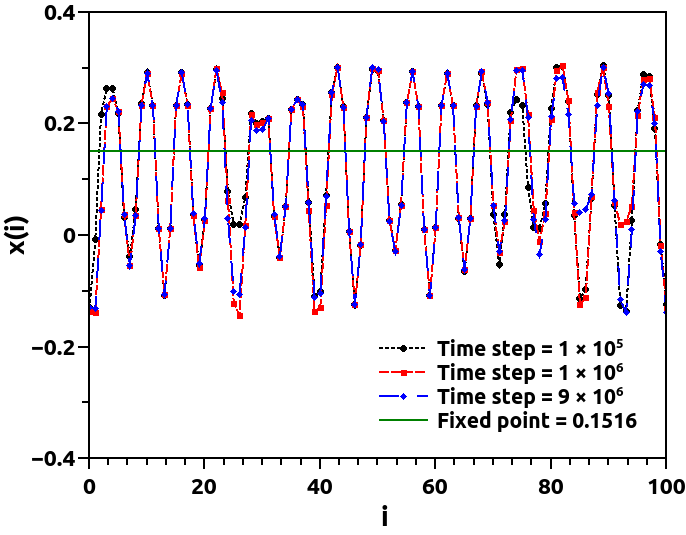}%
\caption{\label{fig:5} Frozen spatial pattern
for $\epsilon=0.4$ and $\beta=-0.69$. 
The plot is done after an even number of time steps.
       The fixed point is also shown as a reference.}
\end{figure}

\begin{figure}
\centering
\includegraphics[width=\linewidth]{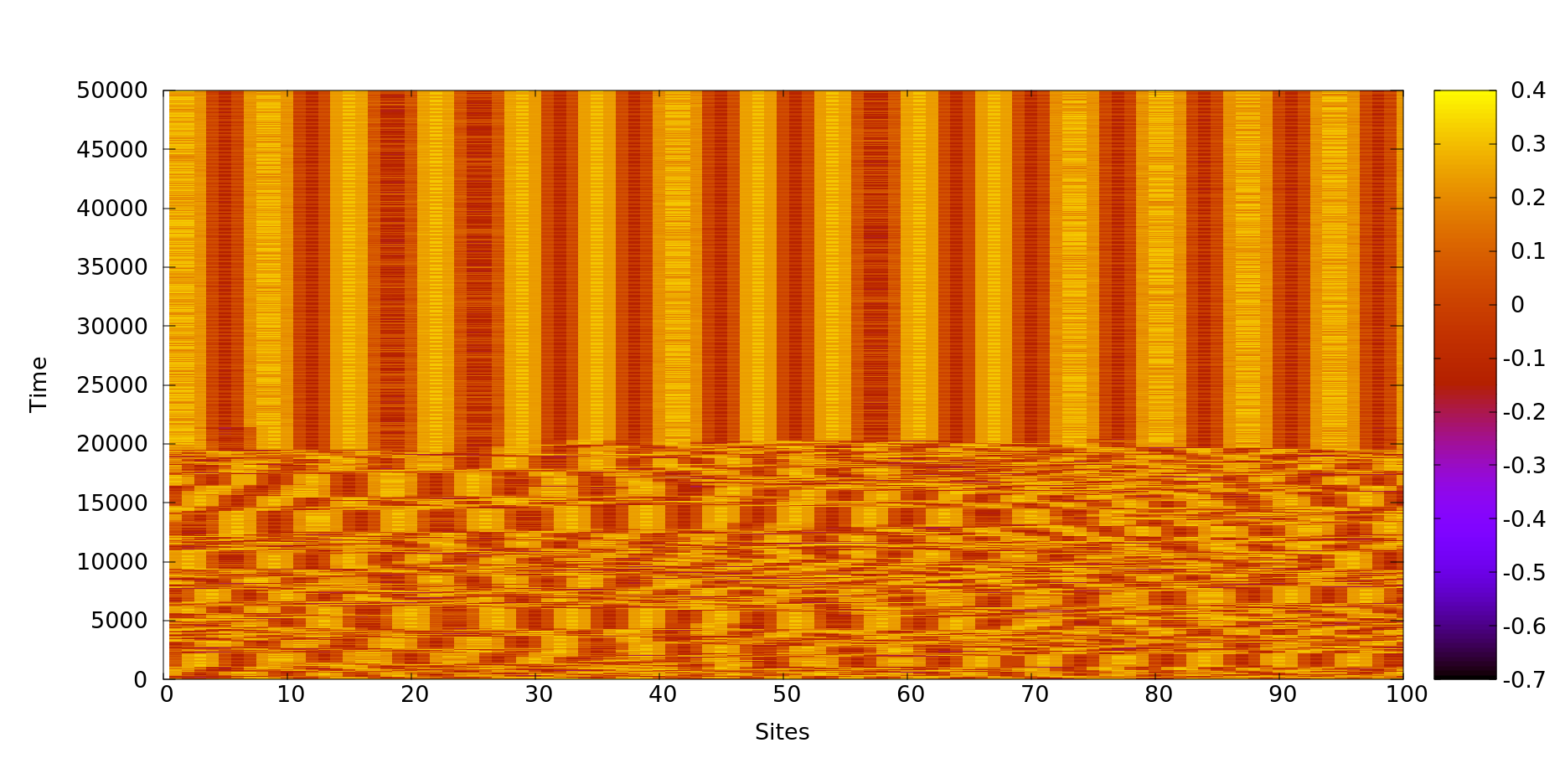}%
	\caption{\label{fig:6} Space-time diagram at critical point $\beta_c=-0.6773$}
\end{figure}

 It is clear that
even though the exact variable value at $i$'th site is not repeated, the sites
with the variable value greater than the largest fixed point $x^*$, 
continue to have a value greater
than $x^*$ and vice versa for $\beta=-0.69$. This behavior
is obtained for $\beta<\beta_c=-0.6773$.

We also plot a space-time diagram at the  point
$\beta=\beta_c=-0.6773$ shown in Fig. 6.
As mentioned above, the
system approaches a state with no spatial period but has a coarse-grained 
period-2
in time. There is no
periodicity in space even in coarse-grained variables. We associate
$s_i(t)=1$ for $x_i(t)>x^*$ and $s_i(t)=-1$ for $x_i(t)<x^*$.
The snapshot of the final spatial pattern 
in these variables is shown in Fig. 7. We observe a predominantly
3-up, 4 -down or 4-up 3-down or 3-up, 3-down pattern. However,
there is no exact periodicity even in coarse-grained sense.
There are infinite such states possible. Clearly, Janssen-Grassberger
conjecture does not apply here. 
The largest fixed point $x^*$  of the map can be found using the bisection method
or other root-finding algorithms.
The system enters a 2-band attractor eventually. So we expect sites
to have the same spin value at all even times and different value at
all odd times. We quantify the transition to a 2-band-attractor
state using two quantifiers, namely
a) Flip rate $F(t)$: The fraction of sites $i$ such that $s_i(2t-2)\ne s_i(2t)$.
b) Persistence $P(t)$: Fraction of sites $i$ such that
$s_i(2t')=s_i(0)$ for all $t'\le t$.

The flip rate $F(t)$ is an indicator of activity
in the lattice at a given time. This is
similar to the density of active sites which is a standard
order parameter for absorbing state in DP class.
We also study
persistence in this system and show that the results are consistent.
We note that if the site $i$ is persistent till time
$T$, it implies that  $s_i(2t)=s_i(0)$ for all $t \leq T$. 
By definition, $P(t)$ decreases monotonically in time. 
It is also possible that it saturates to a finite value
{\it {i.e.}} some sites do not deviate even once from their initial state 
during the entire course of evolution.
Thus nonzero persistence may be due to frozen states in which sites 
do not flip any longer or at least a fraction of sites does not flip
any longer. The transition could be to a fully absorbing state or a partially
absorbing state. The order parameter $F(t)$ helps us
to distinguish between these possibilities.
In this study, both $P(t)$ and $F(t)$ tend
to zero asymptotically at the critical point and a fully absorbing
state is reached. Above the critical 
point, we have a fluctuating state where $F(t) > 0$  
asymptotically and $P(\infty)$ goes to zero since every site flips 
from its initial 
state sooner or later. Below the critical point, $F(\infty)\sim 0$ and some sites
are stuck forever in their initial conditions leading to a nonzero 
asymptotic value
of persistence $P(\infty)>0$.   

We carry out simulations for $N=2\times 10^5$ and average over
at least  $10^3$ 
configurations .
For $\epsilon=0.4$, we observe a clear power-law decay of $F(t)$
as a function of time $t$ at $\beta=\beta_c=-0.6773$. The flip rate 
$F(\infty)$ 
saturates 
for $\beta>\beta_c$ and  $F(\infty)=0$ for $\beta<\beta_c$.
At the critical point, $\beta=\beta_c$,  the asymptotic behavior is expected to
be $F(t) \sim t^{-\delta}$. We indeed observe this behavior in Fig. 8.
The exponent $\delta=0.158$ which is very close to 
directed percolation value of $0.159$ \cite{henkel2008non}.

\begin{figure}
\centering
\includegraphics[width=80mm]{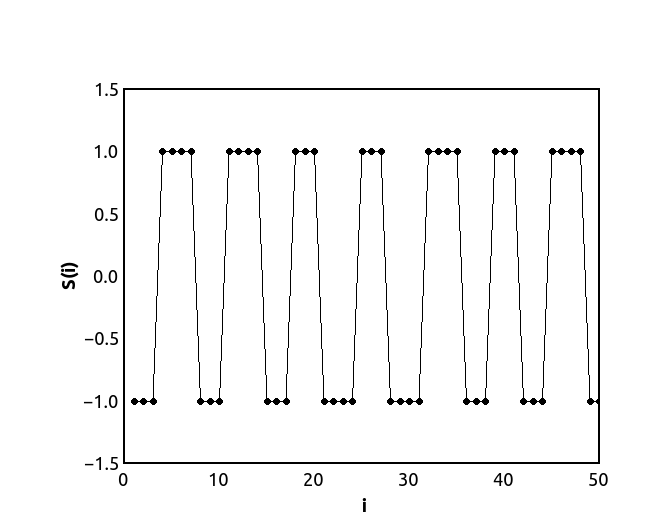}
        \caption{\label{fig:7} 
	Coarse-grained spatial pattern at $\beta = -0.6973$ in the 
	persistent region.}
\end{figure}

\begin{figure}
\centering
\includegraphics[width=80mm]{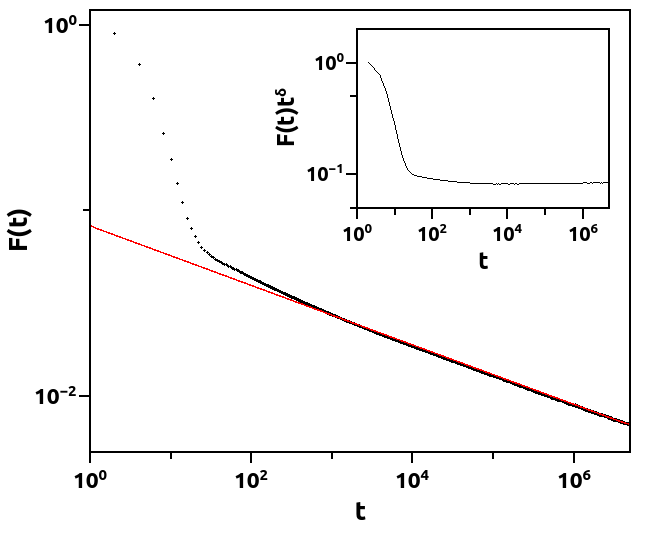}
        \caption{\label{fig:8} Order Parameter $F(t)$ as a function of time $t$
        at the critical point, $\beta_{c} = -0.6773$. 
	We carry out simulation for
        $2 \times 10^{5}$ sites and averaging is carried 
	over $10^3$ configurations.
        Asymptotically, we obtain a clear power-law decay with
	exponent,  $\delta=0.158$. In the inset, $F(t)t^{\delta}$
	is plotted against $t$. Power-law decay is indicated
	from the fact that $F(t)t^{\delta}$  tends to
	a constant.}
\end{figure}

We also find $P(t)$ as a function of time $t$ for $\beta=\beta_c$.
We find $P(t)\sim t^{-\theta}$ with $\theta=1.51$  asymptotically (see Fig. 9). The exponent
$\theta$ is also known as persistence exponent. This value of the exponent
is consistent with persistence exponents obtained
in several other systems showing DP transition except for a couple of models
\cite{matte2016persistence}. Even though the persistence
exponent is not universal, significant sub-class of models showing DP transition
in 1-D have persistence exponent close to $3/2$ \cite{menon2003persistence, hinrichsen1998numerical, 
albano2001numerical, grassberger2009local, fuchs2008local}. 
This is another indicator that the
transition is in DP universality class.

\begin{figure}
\centering
\includegraphics[width=80mm]{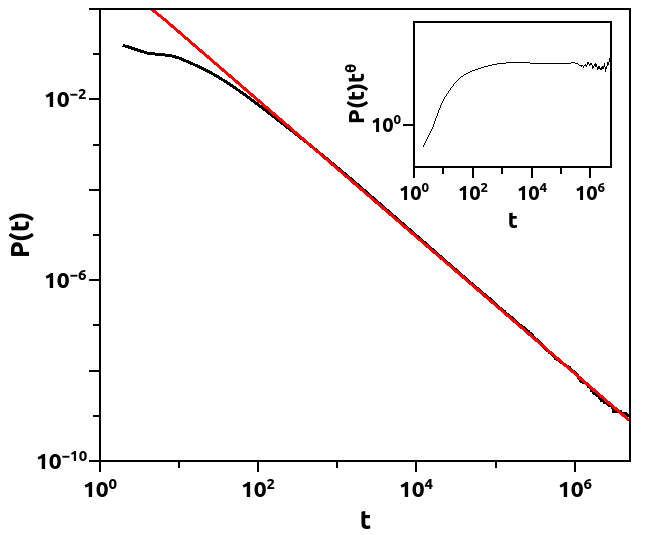}
        \caption{\label{fig:9} Persistence $P(t)$ as a function of time $t$
        at a critical point, $\beta_{c} = -0.6773$. We carry out simulation for
        $2 \times 10^{5}$ sites and average over $5\times 10^{4}$ configurations.
        Asymptotically, we obtain a clear power-law decay with
	persistence exponent,  $\theta=1.51$. In the inset, $P(t)t^{\theta}$
	is plotted against $t$. Power-law decay is indicated
	from the fact that $P(t)t^{\theta}$  tends to
	a constant.}
\end{figure}

We also obtain other exponents such as dynamic exponent $z$ as well as parallel
(temporal) correlation length exponent
$\nu_{\parallel}$ using finite-size scaling and off-critical scaling. We conduct
this exercise for both persistence $P(t)$ as well as for flip rate $F(t)$.

For our order parameter which is flip rate, we expect the
following scaling law to hold, 
\[ F_{N}(t)=t^{-\delta}{\cal F}(t/N^{z}, t \Delta^{\nu_{\parallel}}) \]
where $\Delta=\vert \beta-\beta_c\vert$ is a departure
from the critical point. 
For $N\rightarrow \infty$ $t\rightarrow \infty$, and $\Delta=0$,  ${\cal F}$
tends to a constant and 
 $F(t) \sim t^{-\delta}$.  This fit with $\delta=0.158$
is shown in Fig. 8. 

For persistence we expect that the following asymptotic law to hold:
\[ P_{N}(t)=t^{-\theta}G(t/N^{z}, t \Delta^{\nu_{\parallel}}) \]
where $G$ is the scaling function.
Again in the thermodynamic limit, at the critical point,  
we have $P(t)\sim t^{-\theta}$ and we observe an excellent fit with 
$\theta=1.51$ which matches with estimates of local persistence exponent
in several other one-dimensional models of DP 
\cite{matte2016persistence,hinrichsen1998numerical,fuchs2008local} (See Fig. 9).
For persistence, there is a significant departure at early times.
We fit this departure by incorporating nonlinear correction.
A standard nonlinear correction is given by 
$P(t) \sim C t^{-\theta}(1+c_1 t^{-\gamma} \ldots)$ \cite{barkema2003universality}.
We  ignore
higher-order terms.  The value of $\gamma$ can be found by plotting
$P(t)t^{\theta}$ as a function of  $t^{-\gamma}$ for various values
of $\gamma$. We find good  linear fit for $\gamma=\frac{1}{3}$, 
 $C=13$ and $c_1=-1.7692$. The fit is shown in Fig. 10. 
We retain this correction for  finding values of $z$ and
$\nu_{\parallel}$.

\begin{figure}
\centering
\includegraphics[width=80mm]{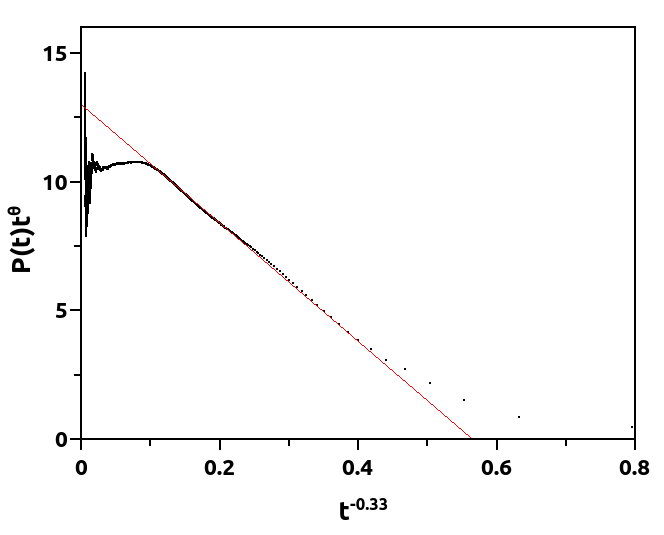}
        \caption{\label{fig:10} $P(t)t^{\theta}$ is plotted against
        $t^{-0.33}$ which gives a good linear fit for the nonlinear correction. }
\end{figure}

For finite-size scaling, we simulate the system 
for various lattice sizes and compute 
the order parameter $F(t)$ as well as $P(t)$ at $\beta=\beta_c$.
The order parameter $F(t)$ is averaged   over $3\times 10^{3}$ configurations,
while  persistence is
averaged over $10^{5}$ configurations. 
The absorbing state is expected to be reached for $t_c=N^z$.
The flip rate scales as 
$F(t_c)=N^{-\delta z}$. We plot $F(t)/F(t_c)$ as 
a function of $t/t_c$ in Fig. 11 and obtain very good scaling collapse. 

\begin{figure}
\centering
\includegraphics[width=80mm]{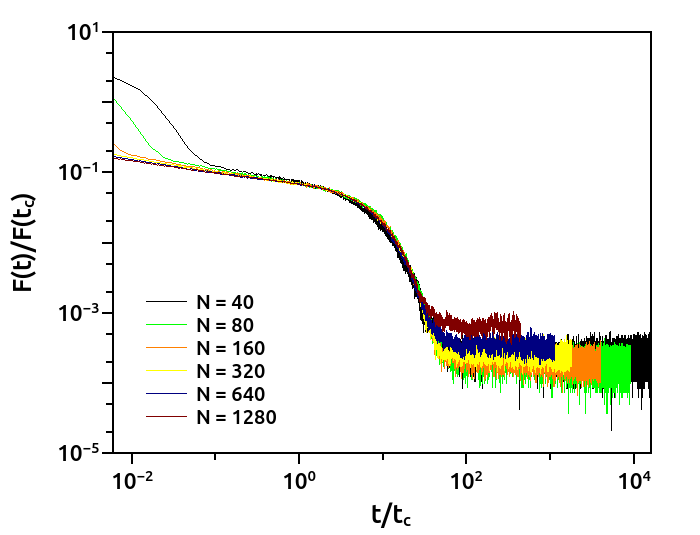}
	\caption{\label{fig:11} Order Parameter $F(t)/F(t_{c})$ as a function of
	$t/t_c$ where $t_c = N^{z}$ for various values of $N$ at critical point $\beta_{c} = -0.6773$. $N$ ranges from 40 to 1280 (from top to bottom).}
\end{figure}

Similarly, we plot $P(t)/P(t_c)$ as a function of $t/t_c$  for various values of $N$ 
in Fig. 12 and obtain excellent scaling collapse.
Thus, we observe an excellent scaling collapse for $z=1.58$ for
the flip rate as well as persistence as shown in Fig. 11 and 12.
These scalings are
consistent with  the expected value of $z$ for 1-D DP class.

\begin{figure}
\centering
\includegraphics[width=80mm]{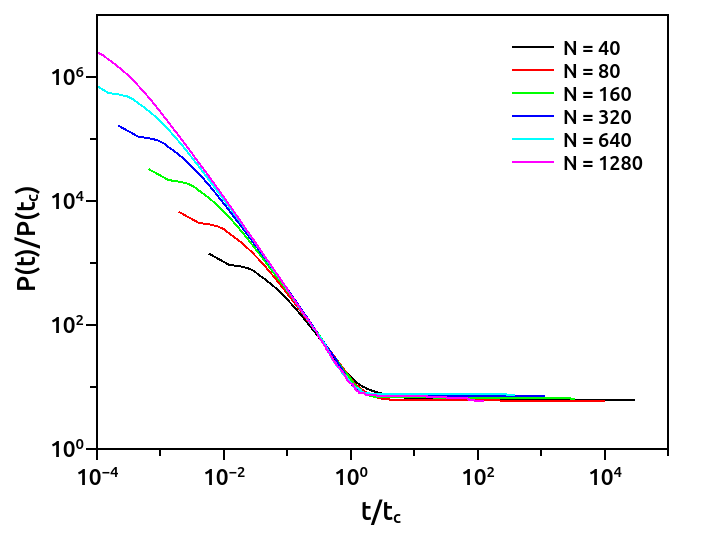}
	\caption{\label{fig:12} Persistence $P(t)/P(t_{c})$ is plotted
	against $t/t_c$ where $t_c = N^{z}$ for various values of $N$ at critical point
	$\beta_{c} = -0.6773$. A clean finite-size scaling is obtained. $N$ ranges from 40 to 1280 (from bottom to top).}
\end{figure}

We also study off-critical scaling behavior for both 
order parameter as
well as persistence to obtain $\nu_{\parallel}$. The size of the lattice is 
large and fixed at
$N=2\times 10^5$ sites. 
Thus finite-size corrections can be neglected.
We carry out extensive averaging  over at least 
$10^3$ configurations. For $\Delta>0$, we average over $10^5$
configurations for persistence.
The critical time scales as $t_c \sim \Delta^{-\nu_\parallel}$.
Thus 
$ F(t_{c}) \sim t_{c}^{-\delta} = \Delta^{\nu_{\parallel}\delta} $.
Plotting $F(t)/F(t_c)$ as a function of $t/t_c$ gives an
excellent scaling collapse for $\nu_{\parallel}=1.73$ (see Fig. 13). This is
an expected value for DP transition \cite{henkel2008non}.

\begin{figure}
\centering
\includegraphics[width=80mm]{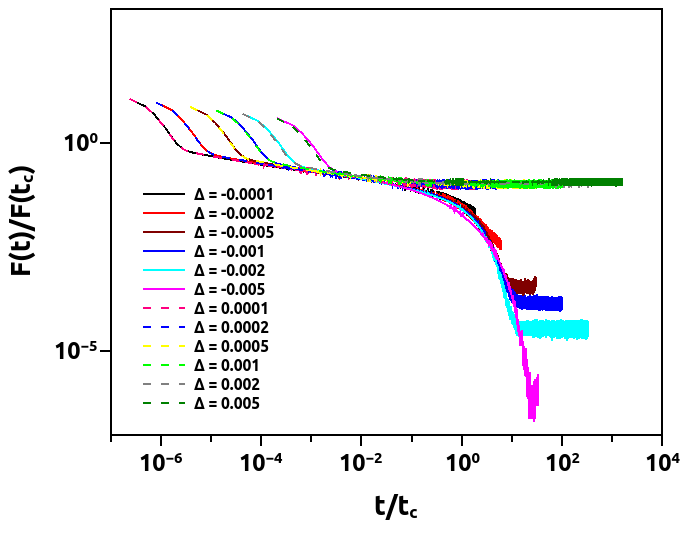}
	\caption{\label{fig:13} Order Parameter $F(t)/F(t_{c})$ is plotted as a 
	function of $t/t_c$ where $t_c = \Delta^{-\nu_{\parallel}}$
        for various values
        of $\Delta=\mid \beta-\beta_{c}\mid$. We carry out off-critical simulation for
        $ 2\times 10^{5}$ sites.}
\end{figure}

We plot $P(t)/P(t_c)$ as a function of $t/t_c$ and obtain an
excellent scaling collapse for $\nu_{\parallel}=1.73$.
The fit is shown in Fig. 14. This value of $\nu_{\parallel}$ is consistent
with the value $1.73$ observed for 1-D DP transitions.

\begin{figure}
\centering
\includegraphics[width=80mm]{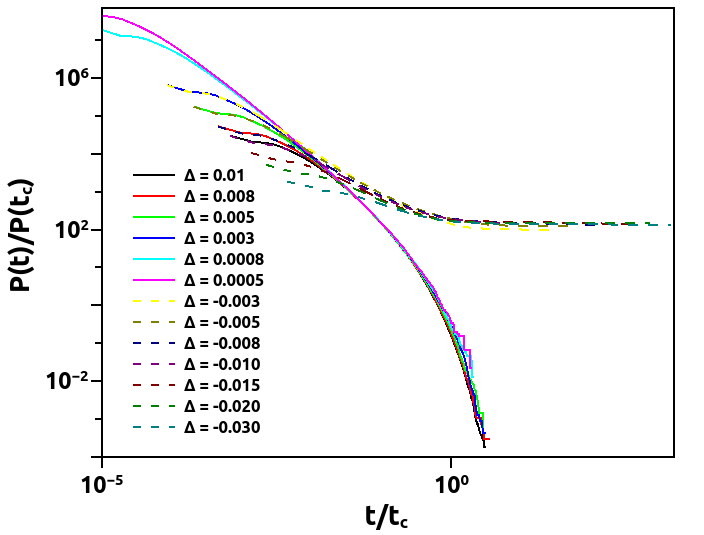}
	\caption{\label{fig:14} Persistence $P(t)/P(t_{c})$ is plotted 
	against $t/t_c$ where $t_c = \Delta^{-\nu_{\parallel}}$
        for various values of
        $\Delta=\mid \beta-\beta_{c}\mid$. We carry out simulations for
        $2 \times 10^{5}$ sites and average over $10^5$ configurations for
        $\Delta > 0$ and $10^3$ for $\Delta<0$.}
\end{figure}

\section{Damage-spreading}
It has been argued that damage-spreading transitions are generically in DP
universality class \cite{grassberger1995damage}. Essentially, we make identical $k$ copies of the
spatially extended system and perturb the central site. We measure the  
difference between all the $N_P=k(k-1)/2$ copies at each instant and sum over those.
We define two quantities, 
$d(t)={\frac{1}{N_P}}
\sum_{l=1}^{k-1}\sum_{m=l+1}^k\sum_{i=1}^n \vert x_i^l(t)-x_i^m(t) \vert
$ and
$D(t)={\frac{1}{N_P}}
\sum_{l=1}^{k-1}\sum_{m=l+1}^k\sum_{i=1}^n \vert s(x_i^l(t))-
s(x_i^m(t)) \vert
$.
We observe that none of these quantities go to zero at our critical
point. Thus the transition is not in  the damage-spreading class.
If we change fraction $p$ of sites in replica, results do not change.

\begin{figure}
	\centering
	\subfloat[]{{\includegraphics[width=8cm]{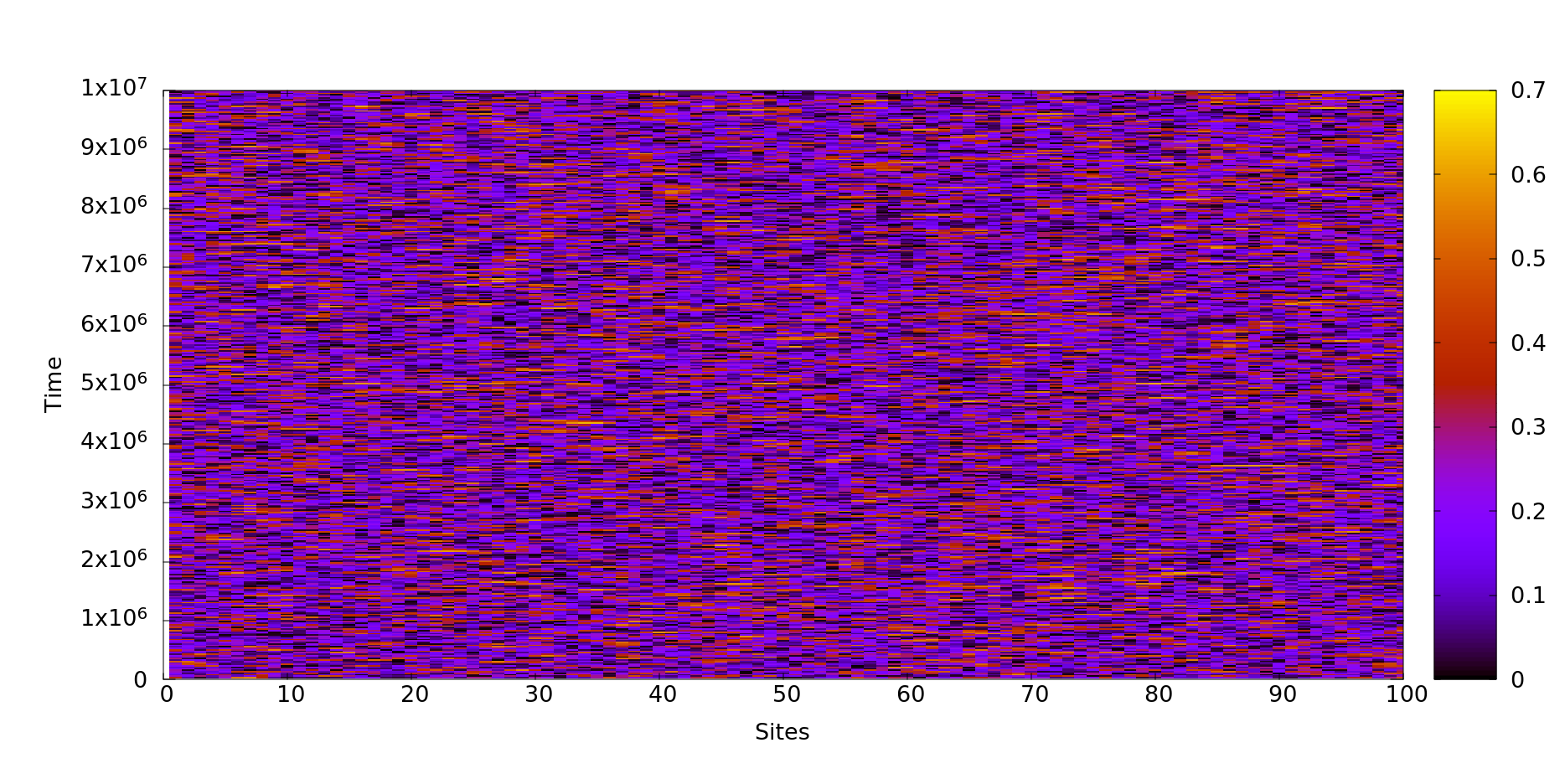}}}%
	\subfloat[]{{\includegraphics[width=8cm]{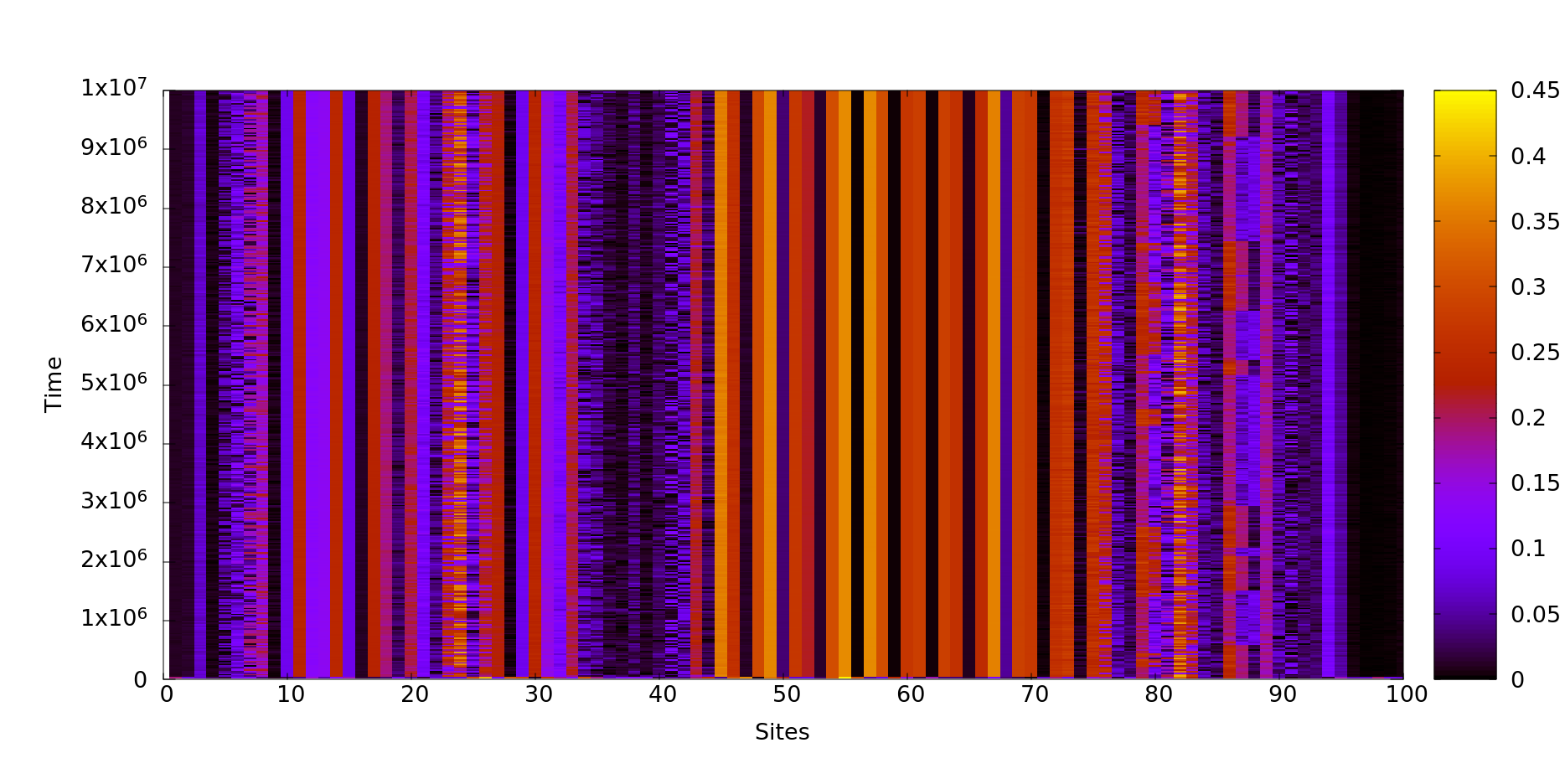}}}%
	\qquad
	\subfloat[]{{\includegraphics[width=8cm]{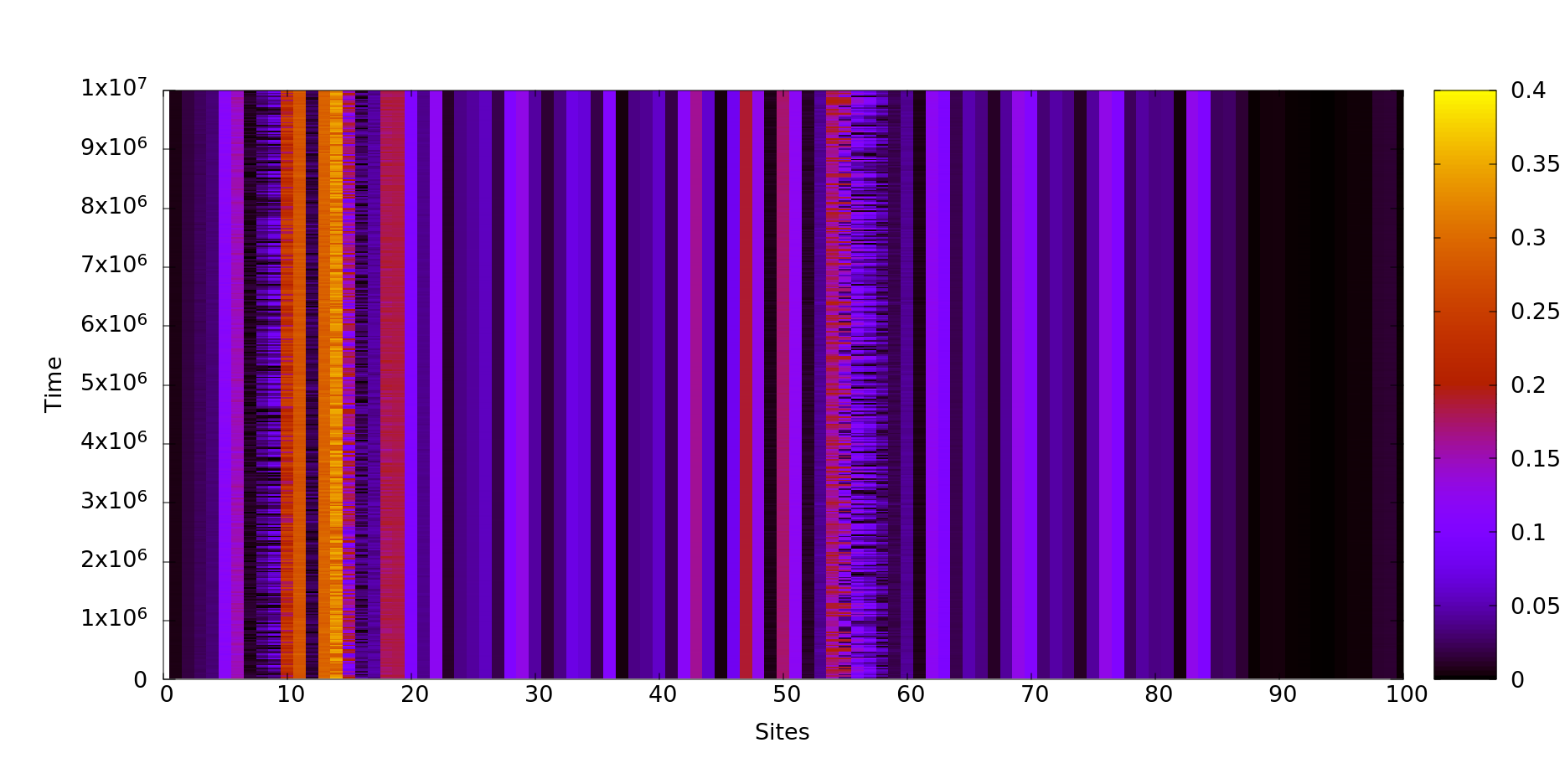}}}%
	\subfloat[]{{\includegraphics[width=8cm]{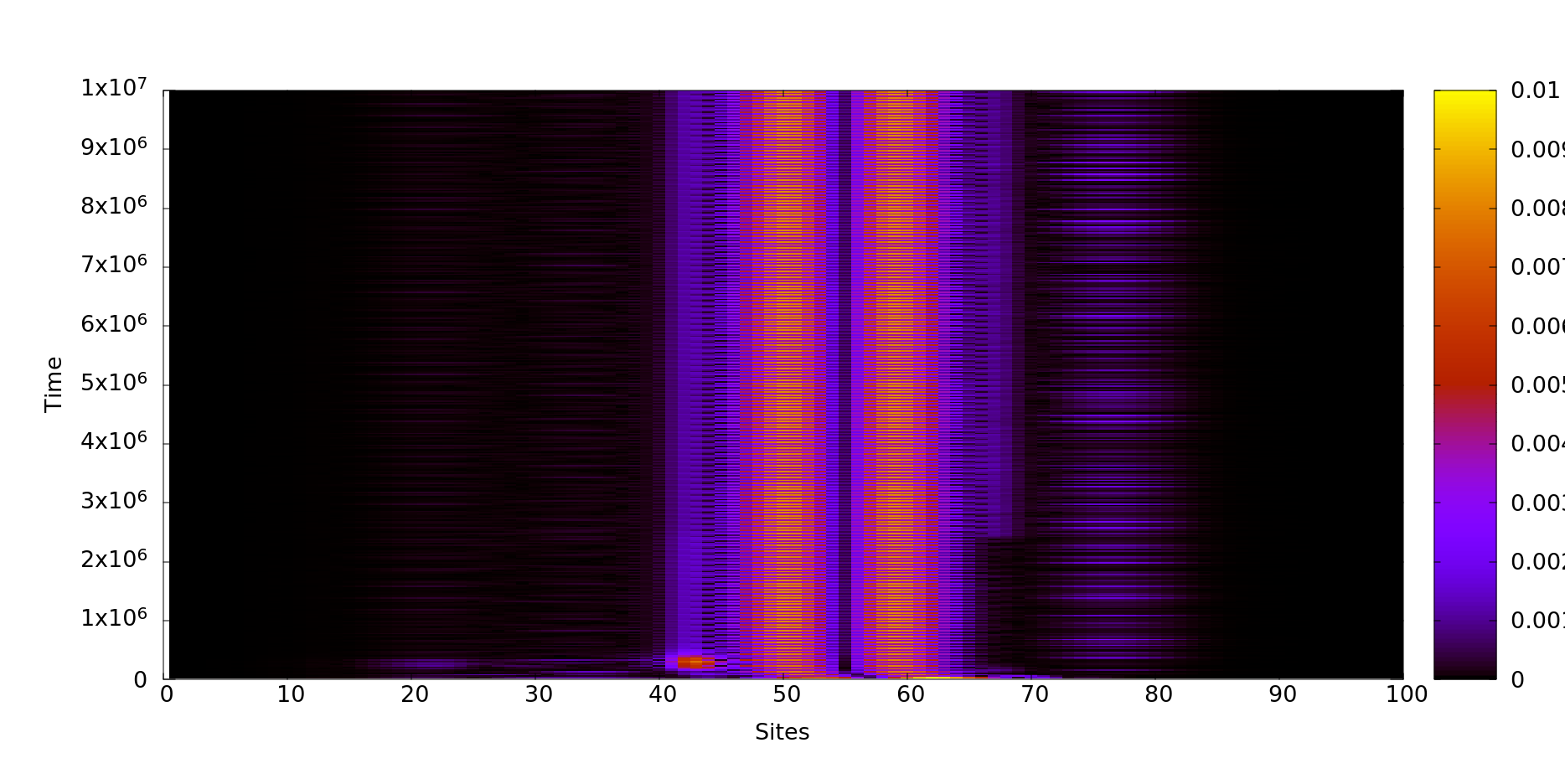}}}%
        \qquad
	\caption{\label{fig:15} Damage spreading at (a) $\beta=-0.63$
	(b) $\beta=-0.6773$ (c) $\beta=-0.73$ (d) $\beta=-0.7771$
	}
\end{figure}

To graphically demonstrate that the
transition is not in 
the damage-spreading class,
we made two identical copies of the system and perturbed the central
site in one of those. Let us denote the variable value at
site i at time $t$ by $x_i(t)$ and the value in its replica by
$y_i(t)$.
Since this is a deterministic system, procedural differences on whether
or not the same set of random numbers is used in simulating both systems
is not relevant.
We plot color-coded difference 
between these two values, {\it{i.e.}}
$\vert x_i(t)-y_i(t) \vert$ as a function of $t$ and $i$ for different
values of $\beta$ (shown in Fig. 15).
The damage spreads almost linearly and spreads to all
sites in the non-persistent region. Even in the persistent region, damage 
does not heal completely. But it does not spread to all sites
and tends to remain localized in the persistent region. 
Thus the transition at $\beta_c$ is not in 
the damage-spreading class if we follow the standard definition of damage
spreading as in \cite{martins1991evidence}. It can be noted that
the damaged sites remain damaged and undamaged sites remain
undamaged below $\beta_c$. Thus the change in the state of damage can be
quantifier which can quantify transition at $\beta_c$. However,
the same information is obtained by flip rate $F(t)$ and no new significant
information is obtained by defining damage in a different manner.
The damage vanishes completely 
only for $\beta<-0.777$ which is far from the critical point.
Both fine-grained damage $d(t)$, as well as coarse-grained
damage $D(t)$, vanish at values of $\beta$ which is much smaller
than the critical point.

We also studied lyapunov exponent. There is no significant change in 
largest lyapunov exponent at $\beta_{c}=-0.6773$. However, the point
at which the damage-spreading transition occurs $(\beta\sim -0.78)$ is the 
one at which the largest lyapunov exponent becomes negative. Thus 
the change in sign of largest lyapunov exponent correlates well with
the damage-spreading transition (see Fig. 16). 

\begin{figure}
\centering
\includegraphics[width=80mm]{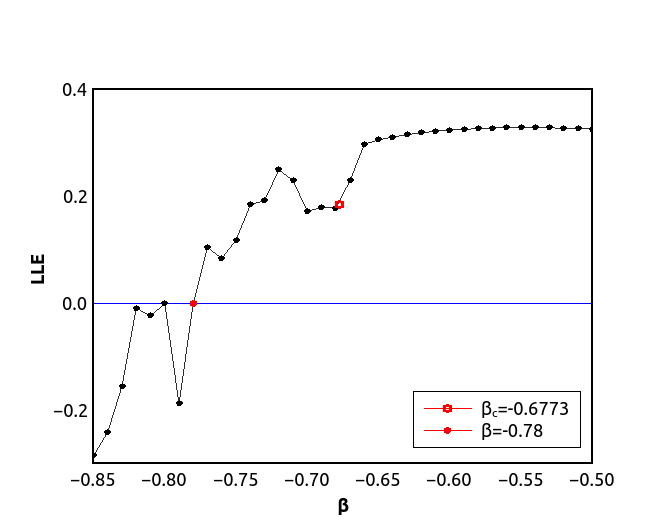}
	\caption{\label{fig:16} Largest lyapunov exponent as a 
	function of $\beta$.
        }
\end{figure}

\section{Discussion}
We study coupled Gauss maps in one dimension.
The system has a transition to an absorbing state
which is a period-2 state in a coarse-grained sense. 
We found no evidence of long-range order
in the system. 
(In fact, for coupled logistic maps in one dimension,
a transition to coarse-grained period-2 state in space and time
is observed.
There is a transition to a state which
shows long-range antiferromagnetic order in space and 
period-2 in time. There are two such absorbing states.
This transition does not belong to DP class but
to a Glauber-Ising class \cite{PhysRevE.87.052905}.)
Usually, a vacuum state is an absorbing state 
for DP transitions and the density of active sites is an
obvious 
order parameter.
Here, we propose the quantity $F(t)$.
This is a fraction of sites which do not return to the same 
band after two time-steps. It shows power-law behavior
at the critical point with  exponent $\delta$ which matches with
systems showing
DP transition. 

Recently, another quantifier known as persistence has
been extremely popular for spotting transitions
to a fully or partially absorbing state. If the flipping rate 
eventually becomes zero, there can be a fraction of sites which 
did not flip even once during evolution. These are known as 
persistent sites. Persistence reaches a finite value asymptotically 
in absorbing state and goes to zero in the active state. At the critical
point, persistence may show power-law decay with exponent known 
as local persistence exponent. In several systems in DP class, this exponent
is found to be $3/2$ or very close to it \cite{hinrichsen1998numerical, 
albano2001numerical, grassberger2009local, fuchs2008local, menon2003persistence}.
We also find local persistence exponent $\theta_l=1.51$
which matches with these models in DP class.

We carry out finite-size scaling as well as off critical scaling to find
other exponents $z$ and $\nu_{\parallel}$.
These exponents are found to be $1.58$ and $1.73$ respectively which put 
the transition firmly in DP universality class. The exponents obtained for 
persistence are consistent with those obtained using order parameter.
We find that the transition is not in  the damage-spreading class.

\section{Acknowledgement}
PMG thanks DST-SERB (Department of Science and Technology, India)  project 
(EMR/2016/006685) 
for financial assistance.




  \bibliographystyle{elsarticle-num} 
  \bibliography{gcml.bib}





\end{document}